\documentclass[11pt]{article}
\usepackage{amsmath, amssymb, amsthm, latexsym, amsfonts, epsfig, color,graphicx}

\begin{document}

\newcommand{\s}{\sigma}
\renewcommand{\k}{\kappa}
\newcommand{\p}{\partial}
\newcommand{\D}{\Delta}
\newcommand{\om}{\omega}
\newcommand{\Om}{\Omega}
\renewcommand{\phi}{\varphi}
\newcommand{\e}{\epsilon}
\renewcommand{\a}{\alpha}
\renewcommand{\b}{\beta}
\newcommand{\N}{{\mathbb N}}
\newcommand{\R}{{\mathbb R}}
   \newcommand{\eps}{\varepsilon}
   \newcommand{\EX}{{\Bbb{E}}}
   \newcommand{\PX}{{\Bbb{P}}}

\newcommand{\cF}{{\cal F}}
\newcommand{\cG}{{\cal G}}
\newcommand{\cD}{{\cal D}}
\newcommand{\cO}{{\cal O}}

\newcommand{\grad}{\nabla}
\newcommand{\n}{\nabla}
\newcommand{\curl}{\nabla \times}
\newcommand{\dive}{\nabla \cdot}

\newcommand{\ddt}{\frac{d}{dt}}
\newcommand{\la}{{\lambda}}

\newtheorem{theorem}{Theorem}
\newtheorem{lemma}{Lemma}
\newtheorem{remark}{Remark}
\newtheorem{example}{Example}
\newtheorem{definition}{Definition}

\title{Predictability in Nonlinear Dynamical Systems
\\with Model Uncertainty }

\author{Jinqiao Duan \\
Department of Applied Mathematics \\
Illinois Institute of Technology \\
Chicago, IL 60616, USA \\ E-mail:  \emph{duan@iit.edu}   }

 \date{July 21, 2008}

\maketitle

\begin{abstract}
Nonlinear systems with model uncertainty are often described by
stochastic differential equations.
 Some   techniques from random dynamical systems are discussed.
 They are relevant to better understanding of solution processes of
  stochastic differential equations and thus may shed lights
  on predictability in nonlinear systems
with model uncertainty.

\medskip
{\bf Key Words:} Stochastic differential equations,  stochastic
parameterizations, predictability, uncertainty, invariant
manifolds, impact of noise
\medskip

%{\bf Mathematics Subject Classifications (2000)}: 34F05, 34C45,
%37H10,  60H10,

\end{abstract}

 \tableofcontents

%%%%%%%%%%%%%%%%%%%%%%%%
%%%%%%%%%%%%%%%%%%%%%%%%

\section{Introduction}

    Nonlinear systems are often influenced by   random fluctuations
such as, uncertainty in specifying initial conditions  or boundary
conditions, external random forcing, and  fluctuating parameters.
In building mathematical models for these nonlinear systems,
sometimes, if not often, less-known, less well-understood, or less
well-observed processes (e.g., highly fluctuating fast or small
scale processes) are ignored due to limitations in our analytical
ability or computational power.

    The limitation of predicting dynamical behavior in nonlinear
systems due to  uncertainty in initial condition has been widely
investigated \cite{GH}. This present article discusses model
uncertainty  in nonlinear systems. This issue has attracted a lot
of attention in geophysical community \cite{Allen, Orrell,
Palmer2006, Farrell, Farrell2, Samelson, Muller,
Saravanan-McWilliams, DelSole, Griffa, Cessi}.

    The uncertainties in simulation may also be regarded as a kind of
model uncertainty. This arises in numerical simulations of
multiscale systems that display a wide range of spatial and
temporal scales, with no clear scale separation. Due to the
limitations of computer power, at present and for the conceivable
future,  not all scales of variability can be explicitly simulated
or resolved. Although these unresolved scales may be very small or
very fast, their long time impact on the resolved simulation may
be delicate (i.e., may be negligible or may have significant
effects, or in other words, uncertain). Thus, to take the effects
of unresolved scales on the resolved scales into account,
representations or parameterizations of these effects   are
required \cite{Berselli}.

Stochastic parameterization of unresolved scales or unresolved
processes leads to stochastic dynamical models in weather and
climate prediction \cite{Palmer2005, Palmer2006, Leith, Imkeller3,
Chorin, Hasselmann, Arn00, Kantz, Stuart, Berloff, DuanNadiga,
Sura3, Williams, Wilks, DuDuan}.

It has been a recent research focus in the dynamical systems
community for better understanding the solution orbits of
stochastic dynamical models \cite{Arnold, Crauel, Oksendal, KS,
DaPrato, Roz, Gar, WaymireDuan}. This is relevant to the issue of
predictability under uncertainty in nonlinear systems, which
concerns about factors and mechanisms for uncertainties of
forecasts and techniques for quantifying and reducing these
uncertainties
 \cite{Palmer2006, Mu2, NorCah81, Boffetta, Smith, Hansen,
Nicolis, Toth, BlomkerDuan}. Various measures have been proposed
in
 quantifying predictability \cite{Kleeman, Majda-Kleeman, Mu1, Schneider},
 and the impact  of measure selection on prediction results
  has also been discussed
 \cite{Orrell2}.

%Mu et al \cite{} have proposed a procedure called conditional
%nonlinear optimal perturbation for quantifying predictability in
%weather and climate dynamics.
%Kleeman  relative entropy

\medskip

We consider the  following   stochastic system defined by
  Ito stochastic differential equations (SDEs) in $\R^n$:
\begin{equation} \label{system}
dX_t   =  b(X_t )d t + \s(X_t )dW(t), \;\; X(0) =  x_0,
\end{equation}
%\begin{equation} \label{system}
% \left \{ \begin{array}{rl}
%dX   =& F(X )d t + B(X )dW(t),\\
%X_0 =& x_0,
%\end{array}
% \right \}
%\end{equation}
where   $b$ and $\s$ are    vector and matrix functions, taking
values in $\R^n$ and $\R^{n\times m}$, respectively. The standard
vector Brownian motion $W(t)$    takes values in $\R^m$. Note that
$n$ and $m$ may be equal or different. We treat $X, X_t, X(t)$ or
$X_t(\om)$ as the same random quantity.

The noise term $\s dW_t$ may be regarded as model uncertainty or
model error. It could be caused by external fluctuations or random
influences, or by a fluctuating coefficients or parameter in the
model. Stochastic parameterization of unresolved scales or
unresolved processes leads to stochastic dynamical systems
\cite{Palmer2005, Palmer, DuanNadiga, Williams, Wilks}. Moreover,
numerical simulation of stochastic partial differential equations
may also lead to SDEs \cite{Millet, Roberts}.

The Brownian motion $W(t)$, or also denoted as $W_t$, is a
Gaussian stochastic process on a underlying probability space
$(\Om, \cF, \PX)$, where $\Om$ is a sample space, $\cF$ is a
$\sigma-$field composed of measurable subsets of $\Om$ (called
``events"), and $\PX$ is a probability (also called probability
measure). Being a Gaussian process, $W_t$ is characterized by its
mean vector (taking to be the zero vector) and its covariance
operator, a $n\times n$ symmetric positive definite matrix (taking
to be the identity matrix). More specifically, $W_t$ satisfies the
following conditions \cite{Oksendal}:

\noindent (a)\ \ \ W(0)=0, \; a.s.\\
(b)\ \ \ W has continuous paths or trajectories, \; a.s.\\
(c)\ \ \ W has independent increments,  \\
(d)\ \ \ W(t)-W(s) $\sim$ $N(0, (t-s)I)$, $t\  and\  s >0 \ and\ t
\geq s \geq 0$, where $I$ is the $n\times n$ identity matrix.

\begin{remark}
(i) The covariance operator here is a constant  $n\times n$
identity matrix $I$, i.e., $Q=I$ and $ Tr(Q)=n$.

(ii) From now on, we   consider two-sided Brownian motion $W_t$,
 $t \in \R$, by means of two independent usual Brownian motions
$W_t^1$ and $W_t^2$ ($t \geq 0$): For $t\geq 0$, $W_t:=W_t^1$,
while for $t<0$, $W_t:=W_{-t}^2$.

 (iii) $W(t)  \sim N(0, |t| I)$, i.e., $W(t)$ has probability
density function $p_t(x) = \frac1{(2\pi t)^{\frac{n}{2}}}
e^{-\frac{x_1^2+...+x_n^2}{2t}}$.

(iv) For every $\a \in (0, \frac12)$, for a.e. $\om \in \Om$,
there exists $C(\om)$ such that

$$ |W(t, \om)-W(s, \om)| \leq
C(\om) |t-s|^{\a},
$$
namely, Brownian paths are H\"older continuous with exponent less
than one half.
\end{remark}

The Euclidean space $\R^n$  has the usual distance
$d(x,y)=\sqrt{\sum_{j=1}^n (x_j-y_j)^2}$, norm   $\|x\|
=\sqrt{\sum_{j=1}^n  x_j ^2}$, and the   scalar product $x\cdot y
= <x, y>= \sum_{j=1}^n x_j y_j$.

This article is organized as follows. After reviewing some basics
about stochastic differential equations in \S \ref{SDE-section},
we discuss random dynamic systems in \S \ref{s2}. Then we consider
the impact of uncertainty and error growth in \S \ref{impact},
residence time, exit probability and predictability in \S
\ref{exit}, and invariant manifolds and predictability in \S
\ref{invariant}. Finally, we discuss nonlinear systems under
non-Gaussian noise and colored noise in \S \ref{nongaussian} and
in \S \ref{colornoise}, respectively.

\section{Stochastic differential equations} \label{SDE-section}

\subsection{Ito and Stratonovich calculus }

Note that the Stratonovich stochastic differential $\s(X) \circ
dW(t)$ and Ito stochastic differential $ \s(X) dW(t)$ are
interpreted through their corresponding definitions of stochastic
integrals \cite{Oksendal}:
$$
\int_0^T \s(X) \circ dW(t) := \mbox{mean-square} \lim_{\Delta t_j
\to 0} \sum_j \s(X(\frac{t_{j+1}+t_j}2)) (W_{t_{j+1}}-W_{t_j}),
$$
$$
\int_0^T \s(X)  dW(t) := \mbox{mean-square} \lim_{\Delta t_j \to
0} \sum_j \s(X(t_j)) (W_{t_{j+1}}-W_{t_j}).
$$
Note the difference in the sums: In Stratonovich integral, the
integrand is evaluated at the midpoint $\frac{t_{j+1}+t_j}2$  of a
subinterval $(t_j, t_{j+1})$, while for Ito integral, the
integrand is evaluated at the left end point $t_j$. See
\cite{Oksendal} for the discussion about the difference in
physical modeling by these two kinds of stochastic differential
equations. There are also dynamical differences for these two type
of stochastic equations, even at linear level \cite{CL}.

  If the integrand $f(t, \om)$ is sufficiently smooth
in time, e.g., H\"older continuous in time in mean-square norm,
with exponent larger than $1$, then both Ito and Stratonovich
integrals coincide; See \cite{Oksendal}, p.39. But in general,
these two Ito and Stratonovich integrals differ. Note that $W_t$
is only H\"older continuous in time \cite{Klebaner}  with exponent
$\a < \frac12$. So that is why the following stochastic integrals
are different:
\begin{eqnarray*}
\int_t^T W_t dW_t =\frac12 (W_t^2-W_T^2)-\frac12 (T-t), \\
 \int_t^T W_t \circ dW_t =\frac12 (W_t^2-W_T^2).
\end{eqnarray*}

Thus we have the two different kinds of SDEs of Ito and
Stratonovich types:
\begin{equation} \label{ito}
dX   =  b(X )d t + \s(X )dW(t), \;\; X(0) =  x_0,
\end{equation}

\begin{equation} \label{strat}
dX   =  b(X )d t + \s(X )\circ dW(t), \;\; X(0) =  x_0,
\end{equation}
 However,
systems of Stratonovich SDEs can be converted to Ito  SDEs and
vice versa  \cite{Kloeden, Oksendal}. In the following we only
consider Ito type of SDEs.

\subsection{Ito's formula  and product rule}

\textbf{Ito's formula in $1-$dimension (scalar case):}

Consider a scalar SDE ($b, \s$ and $W_t$ are all scalars)
\begin{eqnarray*}
 dX_t = b(X_t)dt + \s(X_t) dW_t.
\end{eqnarray*}

Let $g(t, x)$ be a given (deterministic) scalar smooth function.

Ito's formula in differential form is:
\begin{eqnarray}
 d g(t,X_t) = [g_t(t, X_t)+ g_x(t, X_t)b(X_t)
 +\frac12 g_{xx}(t, X_t)\s^2(X_t) ]dt
 + g_x(t, X_t) \s(X_t) dW_t.
\end{eqnarray}
 The term $\frac12 g_{xx}(t, X_t)\s^2(X_t)$ is
called the Ito correction term. Symbolically, we may use the
following rules in manipulating Ito differentials:
\begin{eqnarray*}
 dt dt =dt dW_t =0, \;\;\; dW_t dW_t =dt.
\end{eqnarray*}

 Ito's formula in integral form is:
\begin{eqnarray}
  g(t,X_t) & = & g(0, X_0)
  + \int_0^t[g_t(s, X_s)+ g_x(s, X_s)+\frac12 g_{xx}(s, X_s)\s^2(X_s)
  ]ds    \nonumber \\
& + & \int_0^t g_x(s, X_s) \s(X_s) dW_s.
\end{eqnarray}

The Generator $A$ for this scalar SDE is
$$
A g = g_x b +\frac12 g_{xx} \s^2.
$$

\textbf{Ito's formula in   $n-$dimension (vector case):}

Consider a SDE system
\begin{eqnarray*}
 dX_t = b(X_t)dt + \s(X_t) dW_t,
\end{eqnarray*}
where  $b$ is an $n-$dimensional   vector function, $\s$ is an $n
\times m$ matrix function, and   $W_t(\omega)$ is an
$m-$dimensional Brownian motion.

 Let $g(t, x)$ be a given
(deterministic) scalar smooth function for $x\in \R^n$.

Ito's formula in differential form is:
\begin{eqnarray}
 d g(t,X_t) & = & \{g_t(t, X_t)+ (\grad g(t, X_t))^T b
 +\frac12 Tr[\s \s^T H(g)](t,
 X_t)\} dt \nonumber \\
 & + & (\grad g(t, X_t))^T \s(X_t) dW_t,
\end{eqnarray}
where $ ^T$ denotes transpose matrix,   $H(g)= (g_{x_i x_j})$ is
the $n \times n$ Hessian matrix, and $Tr$ denotes the trace of a
matrix.

Generator $A$  for this SDE system is,
\begin{eqnarray} \label{gtor} A g
= (\grad g)^T b+\frac12 Tr[\s \s^T D^2(g)],
\end{eqnarray}
where the gradient vector of $g$ is $\grad g =(g_{x_1}, \cdots,
g_{x_n})^T$ and the $n\times n$ Hessian matrix of $g$ is
$$
D^2 g := (g_{x_i x_j}).
$$
Symbolically we may also use the rules:
\begin{eqnarray*}
 dt dt =0, \;\; dt dW_t =\textbf{0}, \;\;  dW_t \cdot dW_t =ndt
 =Tr(Q) dt
\end{eqnarray*}
Note that: $Tr (Q) =n$ for $n-$dimensional Brownian motion $W_t$.

Ito's formula in integral form is:
\begin{eqnarray}
  g(t,X_t) & = & g(0, X_0)
  + \int_0^t \{g_t(s, X_s)+ (\grad g(s, X_s))^T b+\frac12 Tr[\s \s^T H(g)](s,
 X_s)\}  ds    \nonumber \\
& + & \int_0^t  (\grad g(s, X_s))^T \s(X_s)   dW_s.
\end{eqnarray}

\begin{remark}
For the above Ito's formula, a somewhat remote connection is the
material derivative of
 the  fluid velocity
\begin{eqnarray}
\ddt g(x, t) = \p_t g + u \cdot \n g,
\end{eqnarray}
where $u=\dot{x}$ is the underlying driving flow.
\end{remark}

\bigskip

\textbf{Stochastic product rule}:

Taking $g =xy$ for a two-dimensional SDE system, we get

\begin{eqnarray}
 d (X_tY_t)= X_t dY_t + (dX_t) Y_t + dX_t dY_t
\end{eqnarray}

\subsection{Estimation of Ito's integrals}

\textbf{Ito isometry}:
\begin{eqnarray}
\EX (\int_0^T f(t,\omega) dW_t)^2  =\EX \int_0^T f^2(t,\omega) dt
\end{eqnarray}

\textbf{Generalized Ito isometry}:

\begin{eqnarray}
 \EX (\int_0^a f(t,\omega) dW_t \int_0^b g(t, \omega)
dW_t)=\EX  \int_0^{a\wedge b} f(t,\omega) g(t,\omega) dt,
\end{eqnarray}
where $a\wedge b  = \min(a, b)$.

\begin{proof}
Look at the case $a=b$.

Denote $I_1 =\int_0^a f(t,\omega) dW_t$ and $I_2 = \int_0^a
g(t,\omega) dW_t$. Note that $I_1 I_2 = \frac12 [(I_1+I_2)^2
-I_1^2-I_2^2]$ and use the isometry property.

For $a \neq b$, say  $a <b$, i.e.,$\min(a, b)=a$. Extend $f$ to
the time interval $[a, b]$ by setting it zero there. Then apply
the above proof.   \qed
\end{proof}

\textbf{Ito isometry in vector case}:

Let $F(t, \om)$ and $G(t, \om)$ be $n\times n$ matrixes, and $W_t$
be $n-$dimensional Brownian motion.

\begin{eqnarray}
\EX  (\int_0^a F(t,\omega) dW_t) \cdot (\int_0^b G(t,\omega) dW_t)
=\EX \int_0^{a \wedge b}  Tr (GF^T)(t,\omega) dt,
\end{eqnarray}
where $\cdot$ denotes the usual scalar product in $\R^n$, $Tr$
denotes the trace of a matrix (i.e. the sum of diagonal entries of
a matrix).

In particular,
\begin{eqnarray}
\EX  \| \int_0^a F(t,\omega) dW_t\|^2 =\EX \int_0^a Tr (F
F^T)(t,\omega) dt ,
\end{eqnarray}

\begin{eqnarray}
 \EX  (\int_0^a F(t,\omega) dW_t) \cdot (\int_0^b F(t,\omega) dW_t)
=\EX \int_0^{a \wedge b}  Tr (F F^T)(t,\omega) dt.
\end{eqnarray}

\textbf{Inequalities involving Ito's integrals}:

By the Ito isometry and the Doob martingale inequality
(\cite{Oksendal}, p.33), we   have, for any constant $\la>0$,

\begin{eqnarray} \label{Martingale-Ito}
\PX (\sup_{t_0\leq t \leq T}|\int_{t_0}^t f(s,\omega) dW_s|\geq
\la ) \leq \frac1{\la^2} \EX \int_{t_0}^T |f (s,\omega)|^2 ds.
\end{eqnarray}

Arnold [1974], p.81:

\begin{eqnarray} \label{Ito-square}
\EX (\sup_{t_0 \leq t \leq T}|\int_0^t f(s,\omega) dW_s|^2 ) \leq
4 \EX \int_{t_0}^T |f(s,\omega)|^2 ds.
\end{eqnarray}

More generally,

\begin{eqnarray} \label{Ito-2k}
\EX  |\int_{t_0}^t f(s,\omega) dW_s|^{2k}   \leq (k(2k-1))^{k-1}
(t-t_0)^{k-1} \EX \int_{t_0}^T |f(s,\omega)|^{2k} ds.
\end{eqnarray}

\subsection{Some examples}

\begin{example}
  \textbf{Langevan equation}

\begin{eqnarray}
dX_t=-b X_tdt + a dW_t,
\end{eqnarray}
where $a, b$ are real parameters, and the initial condition $X_0
\sim \N (0, \s^2)$. The solution is
\begin{eqnarray}
X_t = e^{-bt} X_0 +a e^{-bt} \int_0^t e^{bs} dW_s.
\end{eqnarray}

Note that $\EX X_t =0$ and $X_t$ is a Gaussian process.

\begin{eqnarray}
Cov(X_s,X_t) = \s^2 e^{-b(s+t)}   + \frac{a^2}{2b}[
e^{-b|s-t|}-e^{-b(s+t)}], \\
Cov(X_0,X_t) = \s^2 e^{-bt}   + \frac{a^2}{2b}[
e^{-b|t|}-e^{-bt}], \\
Var(X_t) = \s^2 e^{-2bt}   + \frac{a^2}{2b}[ 1-e^{-2bt}].
\end{eqnarray}

\begin{eqnarray}
Cor(X_s,X_t) =  \frac{Cov(X_s,X_t)}{\sqrt{Var(X_s)}
\sqrt{Var(X_t)}}
\end{eqnarray}

When $\s^2 =\frac{a^2}{2b}$, we have
\begin{eqnarray}
Cov(X_s,X_t) = \s^2  e^{-b|s-t|}, \\
Var(X_t) =   \frac{a^2}{2b}.
\end{eqnarray}
Namely, in this case, $X_t$ is a stationary process.

\end{example}

\begin{example}
  \textbf{Stochastic population model}

  Consider the following  linear scalar SDE with multiplicative noise:
\begin{eqnarray}
dX_t = r X_t dt + \alpha X_t dW_t,
\end{eqnarray}
where $r$ and $\alpha$ are real constants, and $X_t>0, a.s.$.
Rewrite the SDE as
\begin{eqnarray} \label{solve100}
\frac{dX_t}{X_t} = r   dt + \alpha   dW_t.
\end{eqnarray}
Applying  the Ito formula to $\ln X_t$ to obtain
\begin{eqnarray}
d (\ln X_t) =\frac{dX_t}{X_t}- \frac12 \alpha^2 dt.
\end{eqnarray}
That is, $\frac{dX_t}{X_t}= d (\ln X_t) + \frac12 \alpha^2 dt$.
Thus \eqref{solve100} becomes
\begin{eqnarray} \label{solve101}
 d (\ln X_t)   = (r-\frac12 \alpha^2)   dt + \alpha   dW_t.
\end{eqnarray}
Integrating from $0$ to $t$,
\begin{eqnarray*}
 \ln \frac{X_t}{X_0}    = (r-\frac12 \alpha^2) t + \alpha W_t.
\end{eqnarray*}
We hence get the final solution
\begin{eqnarray} \label{solve102}
  X_t = X_0 \exp((r-\frac12 \alpha^2) t + \alpha W_t).
\end{eqnarray}

\end{example}

\begin{example}
\textbf{A linear scalar SDE} \cite{Arnold74, Kloeden}:

\begin{eqnarray}
 dX_t = [a_1(t)X_t +a_2(t)]dt + [ b_1(t)X_t + b_2(t)] dW_t, \;\;
 X_{t_0} \; \mbox{given}.
\end{eqnarray}
The fundamental solution
\begin{eqnarray}
 \Phi_{t,t_0} = \exp[\int_{t_0}^t (a_1(s)-\frac12b_1^2(s))ds
 +\int_{t_0}^tb_1(s) dW_s].
\end{eqnarray}
The general solution
\begin{eqnarray}
  X_t = \Phi_{t,t_0}\; \{ X_{t_0} + \int_{t_0}^t [a_2(s)- b_1(s)b_2(s)]\; \Phi^{-1}_{s, t_0}ds \\
  +  \int_{t_0}^t b_2(s) \; \Phi^{-1}_{s, t_0}dW_s \},
\end{eqnarray}

\end{example}

\begin{example}

\textbf{A    linear system  of SDEs} \cite{Oksendal}:

\begin{eqnarray} \label{nonhom-constant}
 dX_t = [AX_t+f(t)] dt + \sum_{k=1}^m g_k(t) dW_k(t), \;\;
 X_{t_0} \; \mbox{given},
\end{eqnarray}
where $A$ is a constant $n\times n$ matrix, $X(t)$, $f(t)$ and
$g_k(t)$'s are $n-$dimensional vector functions, and $W_k$ are
independent scalar Brownian motions. This is a system with
constant coefficient matrix and additive noise. In this case, we
can find out the solution completely with the help of matrix
exponential.

The fundamental solution matrix for the corresponding linear
system $dX_t =  AX_t dt$  is
\begin{eqnarray} \label{fundamental}
 \Phi_{t, t_0} = e^{A(t-t_0)}.
\end{eqnarray}
The solution for the nonhomogeneous linear system
 with constant coefficient matrix (\ref{nonhom-constant}) is
\begin{eqnarray} \label{soln-constant}
X_t =  e^{A(t-t_0)}\{ X_{t_0}
 + \int_{t_0}^t e^{-A(s-t_0)} f(s)  ds  \\
 + \sum_{k=1}^m \int_{t_0}^t e^{-A(s-t_0)}  g_k(s)  dW_k(s)  \} \\
 = e^{A(t-t_0)} X_{t_0}
 + \int_{t_0}^t e^{A(t-s)} f(s)  ds  \\
 + \sum_{k=1}^m \int_{t_0}^t e^{A(t-s)}  g_k(s)  dW_k(s).
\end{eqnarray}

\end{example}

\begin{example}
\textbf{Stochastic oscillations} \cite{Mao, Oksendal}

\begin{eqnarray}
\ddot{x} + a \dot{x} + b x = \s \dot{W}_t,
\end{eqnarray}
where $a, b, \s$ are real constants, and $W_t$ is a scalar
Brownian motion. This second order SDE may be rewritten as a first
order SDE system:
\begin{eqnarray}
\dot{x} & =& y, \\
 \dot{y} & =& - b x -ay + \s \dot{W}_t.
\end{eqnarray}
In matrix form this becomes
\begin{eqnarray}
\dot{X}   & =& AX + K \dot{W}_t,
\end{eqnarray}
where
\begin{equation*}
A=\left(
    \begin{matrix}
 0 & 1 \\
 -b & -a
    \end{matrix}
\right)
\end{equation*}
and
\begin{equation*}
K=\left(
    \begin{matrix}
 0  \\
 \s
    \end{matrix}
\right).
\end{equation*}
The solution  is
\begin{eqnarray}
 X(t) = e^{At}X(0) + \int_0^t e^{A(t-s)} K dW_s.
\end{eqnarray}

A special case of this model is the stochastic harmonic
oscillator:
\begin{eqnarray}
\ddot{x} +   k x = h \dot{W}_t,
\end{eqnarray}
where $k, h$ are positive constants.
%The solution can be obtained via either of the above two methods.
%(see \cite{Mao}, p.276, for solution via the matrix exponential
%method)
In this case ($a=0$),
\begin{equation*}
A=\left(
    \begin{matrix}
 0 & 1 \\
 -k & 0
    \end{matrix}
\right).
\end{equation*}
Noticing that $A^2=-k I$ with $I$ the $2 \times 2$ identity
matrix, we have
\begin{eqnarray}
e^{At} = \left(
    \begin{matrix}
 \cos (\sqrt{k} t) & \frac1{\sqrt{k}}\sin(\sqrt{k} t)  \\
 -\sqrt{k}\sin(\sqrt{k} t) & \cos (\sqrt{k} t)
    \end{matrix}
\right).
\end{eqnarray}
The final solution for the stochastic harmonic oscillator is
\begin{eqnarray}
 x(t) &=& x_0 \cos (\sqrt{k} t) + \frac{y_0}{\sqrt{k}}\sin(\sqrt{k} t)
 + \frac{h}{\sqrt{k}}\int_0^t \sin(\sqrt{k}(t-s))dW_s ,   \\
 y(t) &=& -x_0\sqrt{k} \sin (\sqrt{k} t) +y_0 \cos(\sqrt{k} t)
 + h \int_0^t \cos(\sqrt{k}(t-s))dW_s .
\end{eqnarray}

\end{example}

%%%%%%%%%%%%%%%%%%%%%%%%%%%%%
%%%%%%%%%%%%%%%%%%%%%%%%%%%%%
%%%%%%%%%%%%%%%%%%%%%%%%%%%%%
\section{Random dynamical systems}\label{s2}

 In this section we introduce some   definitions in stochastic
dynamical systems,  as well as recall some usual notations
  in probability.

%$(\Omega, \mathcal{A}, \mathbb{P})$.

We consider stochastic systems in the state space $\R^n$.
%with the
%usual  metric or distance $d(x,y)=\sqrt{\sum_{j=1}^n
%(x_j-y_j)^2}$, norm or length $\|x\| =\sqrt{\sum_{j=1}^n  x_j
%^2}$, and the usual scalar product $<x,   y>= \sum_{j=1}^n x_j
%y_j$.
All    the sample paths or sample orbits and invariant manifolds
are in this state space.

Some stochastic processes, such as a Brownian motion, can be
described by a canonical (deterministic) dynamical system (see
\cite{Arnold}, Appendix A). A standard Brownian motion (or Wiener
process) $W(t)$ in $\R^n$, with two-sided time $t \in \mathbb{R}$,
is a stochastic process with $W(0)=0$ and stationary independent
increments satisfying $W(t)-W(s) \thicksim \mathcal{N} (0,
|t-s|I)$. Here $I$ is the $n\times n$ identity matrix. The
Brownian motion can be realized in a canonical sample space of
continuous paths passing the origin at time $0$
\[
\Om = C_0(\R, \R^n): =\{\om \in C(\R,\R^n): \om(0)=0 \}.
\]
We identify $W_t(\om)$ with $\om(t)$, namely $W_t(\om) =\om(t)$.
The  convergence concept in this sample space is the uniform
convergence on bounded and closed time intervals, induced by the
following   metric
\[
\rho(\om, \om'):= \sum_{n=1}^{\infty}\frac1{2^n}\;
\frac{\|\om-\om'\|_n}{1+\|\om-\om'\|_n},\; \mbox{where}\;
\|\om-\om'\|_n:=\sup_{-n\leq t \leq n} \|\om(t)-\om'(t)\|.
\]
With this metric, we can define events represented by open balls
in $\Om$. For example, a ball centered at zero with radius $1$ is
$\{\om: \; \rho(\om, 0) < 1 \}$. We define the Borel
$\sigma-$algebra $\mathcal{F}$ as the collection of events
represented by    open balls $A$'s,
 complements of open balls, $A^c$'s, unions and intersections of $A$'s and/or $A^c$'s,
together with the empty event,   the whole event (the sample
 space $\Om$), and all events formed by
 doing the complements, unions and intersections forever in this collection.

 Taking the
(incomplete) Borel $\sigma-$algebra $\mathcal{F}$ on $\Om$,
together with the corresponding Wiener measure $\PX$, we obtain
the canonical probability space $(\Omega, \mathcal{F},
\mathbb{P})$, also called the Wiener space. This is similar to the
game of gambling with a dice, where the canonical sample space is
$\Om_{dice}=\{1, 2, 3, 4, 5, 6  \}$. Moreover, $ \mathbb{E}$
  denotes the mathematical expectation with respect to probability $ \mathbb{P}$.

The canonical \emph{driving} dynamical system describing the
Brownian motion is defined as
\[
\theta(t): \Om \to \Om,\;\;\; \theta(t)\om(s):=\om(t+s)-\om(t),\;
s, t \in \mathbb{R}.
\]
Then $\theta(t)$, also denoted as $\theta_t$, is a homeomorphism
for each $t$ and $(t, \om) \rightarrowtail \theta(t)\om$ is
continuous,   hence measurable. The Wiener measure $\PX$ is
invariant and ergodic under this so-called Wiener shift
$\theta_t$. In summary, $\theta_t$ satisfies the following
properties.
\begin{itemize}
    \item $\theta_0=id$,
    \item $\theta_t\theta_s=\theta_{t+s}$, $\;\;$ for all $s$,
    $t\in\R$,
    \item the map $(t,\omega)\mapsto \theta_t\omega$ is
    measurable and $\theta_t\mathbb{P}=\mathbb{P}$ for all
    $t\in\R$.
\end{itemize}

We now introduce an important concept. A filtration is an
increasing family of information accumulations, called
$\sigma$-algebras, $\mathcal{F}_t$. For each $t$, $\sigma$-algebra
$\mathcal{F}_t$ is a collection of events in sample space $\Om$.
One might observe the Wiener process $W_t$ over time $t$ and   use
  $\mathcal{F}_t$ to represent the information accumulated up to
and including time $t$.  More formally, on  $(\Omega,
\mathcal{F})$, a filtration is a family of $\sigma$-algebras
${\mathcal{F}_s : 0 \leq s \leq t}$ with $\mathcal{F}_s$ contained
in $\mathcal{F}$ for each $s$, and $ \mathcal{F}_s \subset
\mathcal{F}_{\tau}$ for $s \leq \tau$. It is also useful to think
$\mathcal{F}_t$ as the $\sigma$-algebra generated by infinite
union of $\mathcal{F}_s$'s, which is contained in $\mathcal{F}_t$.
So a filtration is often used to represent the change in the set
of events that can be measured, through gain or loss of
information.

For understanding stochastic differential equations from a
dynamical point of view, the natural filtration is defined as a
two-parameter family of  $\sigma$-algebras generated by increments
\[
\mathcal{F}_s^t:= \sigma(\om(\tau_1)-\om(\tau_2): s\leq \tau_1,
\tau_2 \leq t), \;\; s, t \in \mathbb{R}.
\]
This represents the information accumulated from time $s$ up to
and including time $t$. This   two-parameter   filtration allows
us to define forward as well as backward stochastic integrals, and
thus we can solve a stochastic differential equation from an
initial time forward as well as backward in time \cite{Arnold}.

The solution operator for the stochastic system (\ref{system})
with initial condition $x(0)=x_0$ is denoted as $\phi(t, \om,
x_0)$.

The dynamics of the system on the state space $\R^n$, over the
driving flow $\theta_t$ is described by a cocycle. A cocycle
$\phi$ is a mapping:
\[
\phi:\mathbb{R} \times \Omega\times \R^n \to \R^n
\]
which  is
$(\mathcal{B}(\mathbb{R})\otimes\mathcal{F}\otimes\mathcal{B}(\R^n),
\mathcal{B}(\R^n))$-measurable
 such that
\begin{eqnarray*}
&\phi(0,\omega,x)=x \in \R^n,\\
&
\phi(t_1+t_2,\omega,x)=\phi(t_2,\theta_{t_1}\omega,\phi(t_1,\omega,x)),
\end{eqnarray*}
for $t_1,\,t_2\in\mathbb{R},\,\omega\in \Omega,$ and $x\in \R^n$.
Then $\phi$, together with the driving dynamical system, is called
a {\em random dynamical system}. Sometimes we also use $\phi(t,
\om)$ to denote this system.

Under very general smoothness conditions on the drift $b$ and
diffusion $\s$, the stochastic differential system  (\ref{system})
  generates a random dynamical system in
$\R^n$; see \cite{Arnold, Kunita}. Let us see an example.

\begin{example}
Consider a SDE:
$$
dX_{t}=X_{t}dt+dW_{t}, \; X_{0}=x \in \R, \;   t \in \R
$$
The solution is $
X_t(\om)=e^{t}x+\int_{0}^{t}e^{t-\tau}dW_{\tau}(\om)$. Thus the
solution operator is
$$
\varphi(t, \omega, x) :
=e^{t}x+\int_{0}^{t}e^{t-\tau}dW_{\tau}(\om).
$$
Note that
\begin{equation}
\varphi(0,\omega,x)=x.
\end{equation}
Now let us show that
\begin{equation}
\varphi(t+s,\omega,x)=\varphi(t,\theta_ s
\omega,\varphi(s,\omega,x)).
\end{equation}
Indeed, on   one hand,
\begin{eqnarray*}
\varphi(t+s,\omega,x)
=e^{t+s}x+\int_{0}^{t+s}e^{t+s-\tau}dW_{\tau}(\om).
\end{eqnarray*}
On the other hand,
\begin{eqnarray*}
\varphi(t,\theta_s\omega,\varphi(s,\omega,x))=e^{t}\varphi(s,\omega,x)+
\int_{0}^{t}e^{t-\tau}dW_{\tau}(\theta_ s\omega
)=\\\\e^{t}[e^sx+\int_0^se^{s-\tau}dW_\tau]+
\int_0^te^{t-\tau}dW_\tau(\theta_s\omega).
\end{eqnarray*}
Now we only to show the following \emph{Claim}:
$\int_0^te^{t-\tau}dW(\theta_s\omega)
=\int_s^{t+s}e^{t+s-\tau}dW_\tau(\om).$ We prove that both sides
of this claim are identical. In fact, noticing that
$dW_\tau(\theta_s\omega)=d(W_{s+\tau}-W_s)$,
\begin{equation}
\mbox{Left hand side} =
\lim-m.s.\sum_{j}e^{t-\tau_j}(W_{s+\tau_{j+1}}-W_{s+\tau_j})
\end{equation}
\begin{equation}
\mbox{Right hand side} =
\lim-m.s.\sum_{j}e^{t+s-(s+\tau_j)}(W_{s+\tau_{j+1}}-W_{s+\tau_j})
=\lim-m.s. \sum_{j}e^{t-\tau_j}(W_{\tau_{j+1}}-W_{\tau_j})
\end{equation}
Hence the claim is proved.
%But $\{B_t\}$ is a stationary, Gaussian process left = right i.e.
%$B_{s+\tau}-B_s$ has the same distribution as $B_\tau$ (since it
%is Gaussian, mean and variance determine all probability
%properties).  See textbook p18. Problem 2.10
Therefore, the solution operator $\varphi(t,\omega,x)$ satisfies
the cocycle property:
\begin{equation}
\varphi(t+s,\omega,x)=\varphi(t,\theta_ s
\omega,\varphi(s,\omega,x))
\end{equation}

\end{example}

\medskip

We recall some concepts in   dynamical systems.  A \emph{manifold}
$M$  is a set, which locally looks like an Euclidean space.
Namely, a ``patch" of the manifold $M$ looks like a ``patch" in
$\R^n$. For example, curves, torus  and spheres in $\R^3$ are one-
and two-dimensional differentiable manifolds, respectively.
However, a manifold arising from the study of invariant sets for
dynamical systems in $\R^n$, can be very complicated. So we give a
formal definition of manifolds. For more discussions on
differentiable manifolds, see \cite{Marsden, Perko}.

\begin{definition} \textbf{(Differentiable manifold and Lipschitz manifold)}
An n-dimensional differentiable manifold $M$, is a connected
metric space with an open covering $\{U_{\alpha}\}$, i.e,
$M=\bigcup_{\alpha}U_{\alpha}$, such that

(i) for all $\alpha$ , $U_{\alpha}$ is homeomorphic to the open
unit ball in $\R^n$, $B=\{x \in \R^n :\; |x| < 1\}$, i.e., for all
$\alpha$ there exists a homeomorphism of $U_{\alpha}$ onto B,
$h_{\alpha}:U_{\alpha} \rightarrow  B$, and

(ii) if $U_{\a} \cap U_{\b}  \neq \varnothing$ and $ h_{\a}:
U_{\a} \rightarrow B$, $h_{\b}: U_{\b} \rightarrow B$ are
homeomorphisms, then $h_{\a}(U_{\a} \cap U_{\b})$ and
$h_{\b}(U_{\a}\cap U_{\b})$ are subsets of $\R^n$ and the map
\begin{equation} \label{map}
 h=h_{\a} \circ h_{\b}^{-1}: h_{\b}(U_{\a}\cap
U_{\b})  \rightarrow  h_{\a}(U_{\a}\cap U_{\b})
\end{equation}
is differentiable, and for all $x\in h_{\b}(U_{\a}\cap U_{\b})$,
the Jacobian determinant $\det Dh(x)  \neq 0$.

If the map (\ref{map}) is only Lispchitz continuous, then we call
 $M$ an n-dimensional  Lispchitz continuous manifold.
\medskip

Recall that a homeomorphism of A to B is a continuous one-to-one
map of A onto B, $h:A \rightarrow B$,  such that $h^{-1}: B
\rightarrow A$ is continuous.
\end{definition}

\medskip

Just as invariant sets are important building blocks for
deterministic dynamical systems,   invariant sets are basic
geometric objects to help understand stochastic dynamics
\cite{Arnold}. Here we present two different concepts about
invariant sets for stochastic systems: random invariant sets and
almost sure invariant sets.

\medskip

 \begin{definition} \textbf{(Random set)}
 A collection  $M=M(\omega)_{\omega\in\Omega}$, of nonempty closed
 sets $M(\omega)$, $\omega\in\Omega$, contained in $\R^n$, is
 called a random  set if
 \begin{eqnarray*}
 \omega\mapsto\inf_{y\in M(\omega)}d(x, y)
 \end{eqnarray*}
 is a random variable for any $x\in \R^n$.
 \end{definition}

% \begin{definition} \textbf{(Tempered absorbing set)}
% A random set $B(\omega)$ is called an tempered absorbing set of
% $\phi$ if for any bounded set $K \subset  \R^n$ there exists
% $t_K(\omega)$ such that $\forall t\geq t_K(\omega)$
% \begin{eqnarray*}
% \phi\big(t,\theta_{-t}\omega,K\big)\subset B(\omega).
% \end{eqnarray*}
% and for all $\varepsilon>0$
% \begin{eqnarray*}
%  \lim_{t\rightarrow\infty}e^{-\varepsilon
% t}d\big(B(\theta_{-t}\omega)\big)=0, \;\;a.e. \;\omega\in \Omega,
% \end{eqnarray*}
% where $d(B)=\sup_{x\in B}d(x, 0)$, with $0\in \R^n$, is the
% diameter of $B$.
% \end{definition}

\begin{definition} \textbf{(Random invariant set)}
A random set $M(\omega)$ is called an invariant set for a random
dynamical system $\phi$ if
\[
\phi(t,\omega,M(\omega))\subset M(\theta_t\omega), \; \;  t \in
\mathbb{R} \;\; \mbox{and}\;\; \om \in \Om.
\]
\end{definition}

Random stationary orbits \cite{Arnold} and periodic orbits
\cite{Zhao} are special invariant sets.

\begin{definition} \textbf{(Stationary orbit)}
 A random variable $y(\om)$ is called a stationary orbit for
 a random dynamical system $\phi$ if
\begin{eqnarray*}
\phi(t, \om, y(\om)) = y(\theta_t \om),\; a.s.,\;  \mbox{for all
$t$}.
 \end{eqnarray*}
\end{definition}

Let us consider an example.
\begin{example}
  Consider a SDE
\begin{eqnarray} \label{eq000}
du(t)=-u(t)dt+dW(t), u(0)=u_0.
\end{eqnarray}
This SDE     defines a random dynamical system
\begin{eqnarray}\label{eq3}
\phi(t, \om, u_0):=u=e^{-t}u(0)+ \int_0^t e^{-(t-s)}dW(s).
\end{eqnarray}
A stationary orbit of this random dynamical system is given by
\begin{eqnarray}\label{eq4}
Y(\omega)=\int_{-\infty}^0 e^s dW_s(\omega).
\end{eqnarray}
\\
\\
Indeed, it follows from (\ref{eq3}) and (\ref{eq4}) that
\begin{eqnarray}\label{eq5}
\phi(t, \om, Y(\omega))&=&e^{-t}Y(\omega)+\int_0^t e^{-(t-s)}dW_s(\omega)\nonumber\\
&=&e^{-t}\int_{-\infty}^0 e^s dW_s(\omega)+\int_0^t e^{-(t-s)}dW_s(\omega)\nonumber\\
&=&\int_{-\infty}^0 e^{-(t-s)} dW_s(\omega)+\int_0^t e^{-(t-s)}dW_s(\omega)\nonumber\\
&=&\int_{-\infty}^t e^{-(t-s)} dW_s(\omega).
\end{eqnarray}
By (\ref{eq3}) we also see that
\begin{eqnarray}
Y(\theta_t \omega)&=&\int_{-\infty}^0 e^s dW_s(\theta_t\omega)\nonumber\\
&=&\int_{-\infty}^0 e^s dW_{s+t}(\omega)\nonumber\\
&=&\int_{-\infty}^t e^{-(t-s)} dW_s(\omega).
\end{eqnarray}
Thus $\phi(t, \om, Y(\om)) = Y(\theta_t \om)$, i.e.,
$Y(\omega)=\int_{-\infty}^0 e^s dW_s(\omega)$ is a stationary
orbit for the random dynamical system \eqref{eq000}.
 \end{example}

\begin{definition} \textbf{(Periodic orbit)}
 A random process $y(t,\om)$ is called an invariant random
 periodic orbit of period $T$ for
 a random dynamical system $\phi$ if
\begin{eqnarray*}
y(t+T,\om) = y(t, \om), \; a.s.  \\
 \phi(t, \om, y(t_0, \om)) = y(t+t_0, \theta_t
\om),\; a.s.
 \end{eqnarray*}
 for all $t$ and $t_0$.
\end{definition}

\begin{definition} \label{randommanifold} \textbf{(Random invariant manifold)}
If  a random invariant set  $M$ can be represented by a graph of a
Lipschitz mapping
\[
\gamma^\ast(\omega,\cdot): \;  H^+\to H^-, \; \mbox{with direct
sum decomposition} \quad H^+\oplus H^-= \R^n
\]
such that
\begin{eqnarray*}
M(\omega)=\{x^++\gamma^\ast(\omega,x^+),x^+\in H^+\},
\end{eqnarray*}
then $M$ is called a Lipschitz continuous invariant manifold.
\end{definition}

%A solution of the stochastic system (\ref{system}) is denoted as
%$x(t, \om)$ or $x_t(\om)$. A natural phase space for the
%stochastic system (\ref{system}) is $L^2(\Omega, \R^n)$, since at
%each time $t$, solution $x(t, \om)$ is a \emph{point} in this
%space.

%Moreover, we can also view  the stochastic system (\ref{system})
%in the following way: For each sample $\om \in \Omega$, the system
%has phase space (or evolves in) $\R^n$.

%\begin{remark}
%Note that a deterministic variable may be regarded as a random
%variable with the sample space having a single sample:
%$\Omega=\{1\}$. Thus a natural phase space for the deterministic
%system $x' = Ax+ f(x,t)$ is $L^2(\Omega, \R^n)$ but with sample
%space $\Omega =\{ 1\}$. In this case $L^2(\Omega, \R^n)$ can be
%identified with $\R^n$, i.e., identifying each $x \in L^2(\Omega,
%\R^n)$ with $x(1) \in \R^n$.
%\end{remark}

\bigskip

We will also consider   deterministic invariant sets or manifolds,
while the invariance   is  in the sense of almost-sure (a.s.)
\cite{Aubin, Filipovic}.

\begin{definition} \label{almostsure}\textbf{(Almost sure invariant set and manifold)}
A (deterministic) set $M$ in $\mathbb{R}^n$ is called locally
almost surely invariant for (\ref{system}), if for all $(t_0,x_0)
\in \R  \times M $, there exists a continuous local weak solution
$X^{(t_0,x_0)}$ with lifetime $\tau=\tau(t_0,x_0)$, such that
\begin{eqnarray*}
X^{(t_0,x_0)}_{t\wedge \tau} \in M , \;\; \forall  t > t_0, \;\;\;
a.s. \;\;  \om \in \Om,
\end{eqnarray*}
where $t\wedge \tau = \min(t, \tau)$. When $M$ is a manifold, it
is called an almost sure invariant   manifold.
\end{definition}

%\begin{definition}
%A set S is called invariant in almost-sure sense for system
%\ref{system} if , for $(t_0,x_0) \in S$, we have
%\begin{eqnarray*}
%P\{(t,X^{(t_0,x_0)}) \in S\  \forall t > t_0 \}=1
%\end{eqnarray*}
%where $X^{(t_0,x_0)}$ is a solution of \ref{system}, such that
%$X^{(t_0,x_0)}=x_0,t_0 \geq 0 $
%\end{definition}

%%%%%%%%%%%%%%%%%%%%%
%%%%%%%%%%%%%%%%%%%%%%%%
%%%%%%%%%%%%%%%%%%%%%%%%
\section{Impact of model uncertainty and error growth }
\label{impact}

Consider a $n-$dimensional SDE system
\begin{eqnarray}
\label{system2}
 dX_t =  b(X_t)dt + \s(X_t) dW_t,
\end{eqnarray}

A typical application of the Ito's formula for SDEs is to estimate
moments of solutions. For example, for the second moment, by
taking $g = \frac12 \|x\|^2 =\frac12 x\cdot x$.
\begin{eqnarray}
 \frac12 d \|X_t\|^2=d g(X_t)=[X_t \cdot b+ \frac12 Tr (\s
 \s^T)]dt + X_t\s(X_t) dW_t
\end{eqnarray}

Taking mean, we get
\begin{eqnarray} \label{energy}
 \frac12 \frac{d}{dt} \EX \|X_t\|^2
 =  \EX (X_t \cdot b)+ \frac12 \EX \; Tr (\s(X_t) \s^T(X_t))
\end{eqnarray}
This tells us how the fluctuating force affects the evolution of
the   mean energy of the system. The final term $ Tr[\s(X_t)
\s^T(X_t)]$  is the effect of noise on mean energy.

Consider the deterministic system without model uncertainty
\begin{eqnarray}
 dY_t =  b(Y_t)dt,
\end{eqnarray}
Then the solution error $U_t=X_t-Y_t$ satisfies
\begin{eqnarray}
 dU_t = [b(U_t+Y_t)- b(Y_t)]dt+\s(U_t+Y_t) dW_t,
\end{eqnarray}
Thus
\begin{eqnarray} \label{errorgrowth}
 \frac12 \ddt \EX \|U_t\|^2 = \EX (U_t \cdot [b(U_t+Y_t)- b(Y_t)])
  + \frac12 \EX Tr[\s(U_t+Y_t) \s^T(U_t+Y_t)].
\end{eqnarray}
This describes the error growth under uncertainty. The final term
$ Tr[\s(U_t+Y_t) \s^T(U_t+Y_t)]$  is the effect of noise on error
growth.

\medskip

Let us look at an   example.

\begin{example} \textbf{Lorenz system under uncertainty}\\
Consider the Lorenz system with multiplicative noise
\begin{eqnarray*}
dx &=&     (-sx +sy)dt + \sqrt{\eps}\; x   dW_1(t), \\
dy &=&   (rx-y-xz )dt + \sqrt{\eps} \;y   dW_2(t),\\
dz &=&  (-bz+xy )dt + \sqrt{\eps} \; z   dW_3(t),
\end{eqnarray*}
where $W_1$, $W_2$ and $W_3$ are independent scalar Brownian
motions, and $r, s, b, \eps$ are positive parameters. The
classical chaos case is when $r=28, s=10 $ and $b=8/3$.

Let $X :=(x, y, z)^T$. Then by the Ito's formula, we obtain energy
estimate
\begin{eqnarray*}
\frac{1}{2}\frac{d}{dt}\mathbb{E}\|X\|^2 &=& \mathbb{E} [-sx^2
-y^2-bz^2 +(r+s) xy +\frac12 \eps (x^2+y^2+z^2)]  \\
 & \leq & [-\min(s, 1, b)+\frac12 (r+s+\eps)] E\|X\|^2,
\end{eqnarray*}
where we have used the fact that $ xy \leq \frac12 (x^2+y^2) \leq
\frac12 (x^2+y^2+z^2)$. We can see that in this case, the noisy
terms add  ``energy" into the system.

\medskip
Now we consider error growth due to uncertainty. Let
$\hat{X}:=(\hat{x}, \hat{y}, \hat{z})^T$ be the (deterministic)
solution    ($\eps=0$ case),  and let $U=(u, v, w)^T:=X-\hat{X}$
be the error. Then by the Ito's formula, we obtain error growth
estimate
\begin{eqnarray*}
\frac{1}{2}\frac{d}{dt}\mathbb{E}\|U\|^2 &=& \mathbb{E} [-su^2
-v^2-bw^2 +(r+s) uv + \hat{y} uw - \hat{z} uv +\frac12 \eps (u^2+v^2+w^2)]  \\
 & \leq & [-\min(s, 1, b)+\frac12 (r+s+|\hat{y}|+|\hat{z}|+\eps)]
 \; \EX \|U\|^2,
\end{eqnarray*}
where we have used the fact that $\EX(\hat{y} uw)   \leq |\hat{y}|
\;  \EX|uw| \leq  \frac12  |\hat{y}|\; \EX (u^2+w^2)$. Note that
under suitable conditions, this system has a random attractor
\cite{Schmalfuss}.

\end{example}

%%%%%%%%%%%%%%%%%%%%%%%%%%%%%%%%%%%
%%%%%%%%%%%%%%%%%%%%%%%%%%%%%%%%%%%
\section{Residence time, exit probability and predictability}
\label{exit}

We start with a SDE system
\begin{eqnarray}  \label{sdes}
 dX_t = b(X_t)dt + \s(X_t) dW_t,\;\;\; X_0 \; \mbox{given}
\end{eqnarray}
where  $b$ is an $n-$dimensional   vector function, $\s$ is an $n
\times m$ matrix function, and   $W_t(\omega)$ is an
$m-$dimensional Brownian motion. The generator for this SDE is a
linear second order differential operator as in \S
\ref{SDE-section}
\begin{eqnarray}  \label{generator}
A g = (\grad g)^T b+\frac12 Tr[\s \s^T D^2(g)],
\end{eqnarray}
where $D^2$ is the Hessain differential matrix and $Tr$ denotes
the trace.

For a bounded domain $D$ in $\R^n$, we can consider the exit
problem of random solution trajectories of \eqref{sdes} from $D$.
To this end, let $\p D$ denote the boundary of $D$ and let
$\Gamma$ be a part of the boundary $\p D$. The escape probability
$p(x,y)$ is the probability that the trajectory of a particle
starting at $(x,y)$ in $D$ first hits $\p D$ (or escapes from $D$)
at some point in $\Gamma$, and $p(x,y)$ is known to satisfy
(\cite{Lin, Schuss, BrannanDuanErvin} and references therein)
\begin{eqnarray}
  A p  & = &  0,  \label{eqn3} \\
    p|_{\Gamma} & = & 1,         \label{eqn4}  \\
    p|_{\p D - \Gamma} & = & 0. \label{eqn5}
\end{eqnarray}
Suppose that initial conditions (or initial particles) are
uniformly distributed over $D$. The average escape probability $ P
$ that a trajectory will leave $D$ along the subboundary $\Gamma$,
before leaving the rest of the boundary, is given by (e.g.,
\cite{Lin, Schuss})
\begin{equation}
      P =  \frac{1}{|D|} \int\int_D p(x,y) dxdy,
      \label{average}
\end{equation}
where $|D|$ is the area of domain $D$.

The residence time of a particle initially at $(x,y)$ inside $D$
is the time until the particle first hits $\p D$ (or escapes from
$D$). The mean residence time $u(x,y)$ is given by (e.g.,
\cite{Schuss, Naeh, Risken} and references therein)
\begin{eqnarray}
    A u  & = &  -1,   \label{eqn6}\\
    u|_{\p D}  & = &    0.  \label{eqn7}
\end{eqnarray}

\emph{Relevance to predictability problem.} For low dimensional
SDE systems, such as the Lagrangian dynamical model for fluid
particles in random fluid flows or other truncated model like the
Lorenz model, the exit probability and mean residence time may be
computed by deterministic partial differential equations solvers
\cite{BrannanDuanErvin}. Be selecting the above domain $D$
appropriately, say corresponding to observational data (``data
domain"), we may determine predictability time window, by
monitoring when the system exits the data domain.

%%%%%%%%%%%%%%%%%%%%%%%%%%%%%%%%%%%
%%%%%%%%%%%%%%%%%%%%%%%%%%%%%%%%%%%
\section{Invariant manifolds  and predictability} \label{invariant}

Invariant manifolds provide geometric structures  that describe
  dynamical behavior of nonlinear     systems. Dynamical
  reductions to attracting invariant manifolds  or dynamical
  restrictions to other (not necessarily
  attracting) invariant manifolds are often sought to gain
  understanding of nonlinear dynamics.

There have been recent  works on invariant manifolds for
stochastic  differential equations \cite{Arnold, Wanner2,
DuanLuSchm, DuanLuSchm2}. Random invariant manifolds in the sense
of Definition \ref{randommanifold} are difficult to obtain, even
locally in state space. But almost sure invariant manifolds in the
sense of Definition \ref{almostsure} may be determined, locally in
state space (which also means for finite time in evolution), for
some SDE systems, by a method  of solving first order
deterministic partial differential equations \cite{DuDuan1}.

We consider the  following   stochastic system defined by
  Ito stochastic differential equations in $\R^n$:
\begin{equation} \label{system3}
dX   =  b(X )d t + \s(X )dW(t), \;\; X(0) =  x_0,
\end{equation}
%\begin{equation} \label{system}
% \left \{ \begin{array}{rl}
%dX   =& F(X )d t + B(X )dW(t),\\
%X_0 =& x_0,
%\end{array}
% \right \}
%\end{equation}
where again $b$ and $\s$ are    vector and matrix functions in
$\R^n$ and $\R^{n\times n}$, respectively,  and $W(t)$ are
standard vector Brownian motion  in $\R^n$.  We also assume that
  $b( \cdot) \in C^1(\R^n;\R^n) $ and
  $\s( \cdot) \in C^1(\R^n;\R^{n \times n})$.

For the  nonlinear stochastic system (\ref{system3}), we study
deterministic almost sure  invariant manifolds,   which are not
necessarily   attracting. We reformulate the local invariance
condition as invariance equations, i.e., first order partial
differential equations, and then solve these equations by the
method of characteristics. Although the local invariant manifold
is deterministic, the restriction of the original stochastic
system on this deterministic local invariant manifold is still a
stochastic system but with reduced dimension.

\bigskip

  We are going to derive representations of   invariant finite dimensional
manifolds in terms of $b$ and $\s $,  by using the tangency
conditions for a deterministic $C^2$ smooth manifold  (a
supersurface) $M$ in $\R^n$:
\begin{eqnarray}
\mu(\om,x) := b(\om, x) &-& {1 \over
2}\sum_j[D \s^j(\om, x)]\; \s^j(\om, x) \in T_xM,    \label{c1}    \\
\s^j(\om, x) &\in& T_xM , \;\;    j=1, \cdots, n,  \label{c2}
\end{eqnarray}
where $D$ represents Jacobian operator and $\s^j$ is the $j-$th
column of the matrix $\s$. The above   tangency conditions are
shown to be equivalent to   almost sure   local invariance of
manifold $M$; see e.g., \cite{Filipovic, Aubin}.

The almost sure invariance conditions (\ref{c1})-(\ref{c2}) for
manifold $M$ mean that the $n+1$ vectors, $\mu$ and $\s^j, \;
j=1, \cdots, n$, are tangent vectors to $M$. Namely, these $n+1$
vectors are orthogonal to the normal vectors of manifold $M$.

In other words, if the normal vector  for $M$ at $x$ is $N(x)$,
then the almost sure invariance conditions (\ref{c1})-(\ref{c2})
become the following \emph{invariance equations} for manifold $M$:
For all $x\in M$,
\begin{eqnarray}
 <\mu (x), N(x)> & =& 0,    \label{c3}     \\
<\s^j(x),  N(x)>  &= & 0, \;\;    j=1, \cdots, n,  \label{c4}
\end{eqnarray}
where, as before, $<\cdot, \cdot>$ denotes the usual scalar
product in $\mathbb{R}^n$.

Invariant manifolds are usually represented as graphs of some
functions in $\R^n$. By investigating the above invariance
equations (\ref{c3})-(\ref{c4}), we may be able to find some local
invariant manifolds $M$ for the stochastic system (\ref{system3}).

The goal for this section is to present a method to find some of
these local invariant manifolds. Although  the following result
and example are stated for a codimension $1$ local invariant
manifold, the idea extends to other lower dimensional local
invariant manifolds, as long as the normal vectors $N(x)$ (or
tangent vectors) may be represented; see tangency conditions
(\ref{c3})-(\ref{c4}) above and (\ref{c5})-(\ref{c6}) below.

\medskip

\noindent \textbf{Local almost sure invariant manifold}:  \\
 Let the local invariant manifold $M$ for the stochastic dynamical system
 (\ref{system3}) be represented as a
graph defined by the algebraic equation
\begin{eqnarray} \label{G}
M: \;\;\; G(x_1 , \cdots , x_n) = 0.
\end{eqnarray}
Then $G$ satisfies a system of first order (deterministic) partial
 differential equations and the local invariant manifold $M$ may
 be found by solving these partial
 differential equations by the method of characteristics. By
 restricting the original dynamical system (\ref{system3}) on this
 local invariant manifold $M$, we obtain a locally valid, reduced lower
 dimensional system.

\bigskip

In fact, the normal vector to this graph or surface is, in terms
of partial derivatives, $\nabla G(x) =( G_{x_1} , \cdots , G_{x_n}
) $. Thus the invariance equations (\ref{c3})-(\ref{c4}) are now
\begin{eqnarray}
 <\mu (x), \nabla G(x)> & =& 0,    \label{c5}     \\
 <\s^j(x), \nabla G(x)> &= & 0, \;\;    j=1, \cdots, n,  \label{c6}
\end{eqnarray}
This is a system of first order partial differential equations in
$G$. We apply the method of characteristics to solve for $G$, and
therefore obtain the invariant manifold $M$, represented by a
graph in state space $\mathbb{R}^n$:
 $G(x_1 , \cdots , x_n) = 0$.\\

\noindent {\em Method of Characteristics.}  Consider a first order
partial differential equation for the unknown scalar function $u$
of n variables $x_1 , \cdots , x_n$
\begin{eqnarray} \label{first}
\sum_{j=1}^n a_i(x_1, \cdots, x_n) u_{x_i} = c(x_1, \cdots, x_n),
\end{eqnarray}
with continuous coefficients $ a_i$'s and $c$.

Note that the solution  surface $u=u(x_1,..,x_n,t)$ in
$x_1...\cdots x_nu-$space has normal vectors $N: =(u_{x_1},
\cdots, u_{x_n}, -1)$. This partial differential equation implies
that the   vector $ V= :(a_1, \cdots, a_n, c)$  is perpendicular
to this normal vector and hence must lie in the tangent plane to
the graph of $z=u(x_1, \cdots, x_n)$.

In other words, $ (a_1, \cdots, a_n, c)$  defines a vector field
in $\mathbb{R}^n$, to which graphs of the solutions must be
tangent at each point \cite{McOwen}. Surfaces that are tangent at
each point to a vector field in $\mathbb{R}^n$ are called
\emph{integral surfaces }  of the vector field. Thus to find a
solution of equation (\ref{first}), we should try to find integral
surfaces.

How can we construct integral surfaces? We can try using the
characteristics curves that are the integral curves of the vector
field. That is, $X=(x_1(t), \cdots, x_n(t) )$ is a
\emph{characteristic} if it satisfies the following system of
ordinary differential equations:
\begin{eqnarray*}
{dx_1 \over dt} & = & a_1(x_1,\cdots,x_n), \\
& \cdots &\\
{dx_n \over dt}& = & a_n(x_1,..,x_n),   \\
{du \over dt} & = & c(x_1,..,x_n).
\end{eqnarray*}
A smooth union of characteristic curves is an integral surface.
There may be many integral surfaces. Usually an integral surface
is determined by requiring it to contain (or pass through) a given
initial curve  or an $n-1$ dimensional manifold $\Gamma$:
\begin{eqnarray*}
x_i  & = & f_i(s_1,..,s_{n-1}),i=1..n\\
u  & = & h(s_1,..,s_{n-1})
\end{eqnarray*}
This generates an $n$-dimensional integral manifold $M$
parameterized by $(s_1,..,s_{n-1}, t)$. The solution $u(x_1,
\cdots, x_n)$ is obtained by solving for $(s_1,..,s_{n-1},t)$ in
terms of variables $(x_1, \cdots, x_n)$.

%More generally,if we are given initial conditions on an
%k-dimensional manifold $S$ ($k \leq n$):
%\begin{eqnarray*}
%x_i  & = & f_i(s_1,..,s_k),i=1..n\\
%z  & = & h(s_1,..,s_k)
%\end{eqnarray*}
%This generates an (k+1)-dimensional integral manifold M
%parameterized by $(s_1,..,s_k,\tau)$. The solution
%$u(x_1,..,x_n,t)$ is obtained by solving for $(s_1,..,s_k,\tau)$
%in terms of variables $(x_1,..,x_{k+1})$.
%\begin{eqnarray*}
%s_i  & = & g_i(x_1,..,x_{k+1}),i=1..k\\
%\tau  & = & g(x_1,..,x_{k+})
%\end{eqnarray*}
%Then we plug the above results into the initial condition, we get
%:
%\begin{eqnarray*}
%x_i & = & x_i(x_1,..,x_{k+1}),i=k+2..n
%\end{eqnarray*}
%So the (k+1)-dimensional integral manifold M is represented as
%\begin{eqnarray*}
%z=h(g_1(x_1,..,x_{k+1}),..,g_k(x_1,..,x_{k+1}))=0\  or\\
%G(x_1,..,x_{k+1},x_{k+2}(x_1,..,x_{k+1}),..,x_n(x_1,..,x_{k+1}))=0.
%\end{eqnarray*}

\begin{remark}  If initial data $\Gamma$ is non-characteristic,
i.e., it is nowhere tangent to the vector field $V=(a_1, \cdots,
a_n, c)$, and $a_1, \cdots, a_n, c$ are $C^1$ (and thus locally
Lipschitz continuous), then there exists   a unique integral
surface $u=u(x_1, \cdots, x_n)$ containing $\Gamma$, defined at
least locally near $\Gamma$.
\end{remark}

\bigskip

Now applying the above method of characteristics to
(\ref{c5})-(\ref{c6}), we obtain a solution $G=G(x_1, \cdots,
x_n)$.  However, the local invariant manifold  $M$ that we look
for is represented by the equation
$$
G(x_1, \cdots, x_n) =0.
$$
Therefore, a skill is needed to make sure that the solution
$G=G(x_1, \cdots, x_n)$ actually penetrates the plane $G=0$ in the
$x_1\cdots x_n G-$space; see Fig. 1. This needs to be achieved by
selecting appropriate initial data $\Gamma$.  The invariant
manifold $M$ we thus obtain is defined at least locally near the
initial data $\Gamma$.

\begin{figure}
\psfig{figure=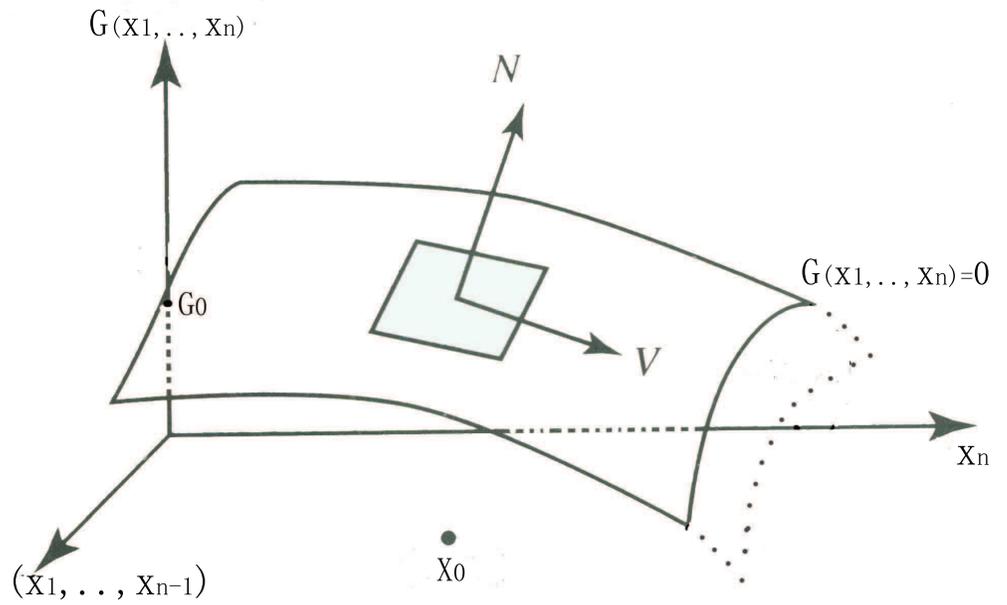,width=6.2in,bbllx=1bp,bblly=1bp,bburx=280bp,bbury=280bp}
    \caption{Local invariant manifold $M$ is represented
     by the equation $G(x_1, \cdots, x_n)=0$
   in the $x_1 \cdots x_n-$ space. Namely, $M$ is the intersection of
    the  surface $G=G(x_1, \cdots, x_n)$ with the  plane $G=0$
    in $x_1 \cdots x_nG-$space. Here
    $G(x_1, \cdots, x_n)$ is the solution of (\ref{c5})-(\ref{c6}) via the
    method of characteristics. Note that $N  =(u_{x_1},
\cdots, u_{x_n}, -1)$ and  $ V=  (a_1, \cdots, a_n, c)$. }
\end{figure}

\bigskip

\emph{Relevance to predictability problem.} When a SDE system
starts to evolve inside a local almost sure invariant manifold
$M$, it remains inside the manifold for a certain time period $0<
t <T$. As determined above, this manifold holds solutions for the
system, the time period $T$ may be taken as a lower bound of the
predictability time scale.

%%%%%%%%%%%%%%%%%%%%%%%%%%%%%%%%%%%
%%%%%%%%%%%%%%%%%%%%%%%%%%%%%%%%%%%
\section{Systems driven by non-Gaussian noise  } \label{nongaussian}

Although Gaussian processes like Brownian motion have been widely
used in modeling fluctuations in geophysical modeling, it turns
out that many physical phenomena  involve with non-Gaussian Levy
motions. For instance, it has been argued that diffusion by
geophysical turbulence \cite{Shlesinger} corresponds, loosely
speaking, to a series of  ``pauses", when the particle is trapped
by a coherent structure, and ``flights" or ``jumps" or other
extreme events, when the particle moves in the jet flow.
Paleoclimatic data \cite{Dit} also indicates such irregular
processes.

Levy motions are thought to be appropriate models for non-Gaussian
processes with jumps \cite{Sato-99}. Let us recall that  a
L\'evy motion $L(t)$ has  independent and stationary increments,
i.e., increments $\Delta L (t, \Delta t)= L(t + \Delta t)-  L(t)$
are stationary (therefore $\Delta L$ has no statistical dependence
on $t$) and independent for any non overlapping time lags $\Delta
t$. Moreover, its sample paths are only continuous in probability,
namely, $\PX (|L(t)-L(t_0)| \geq \delta) \to 0$ as $t\to t_0$ for
any positive $\delta$. This continuity is weaker than the usual
continuity in time.

This generalizes the Brownian motion $B(t)$, as $B(t)$ satisfies
all these three conditions.  But \emph{Additionally},  (i) Almost
every sample path of the Brownian motion     is  continuous in
time in the usual sense and (ii) Brownian motion's increments are
Gaussian distributed.

SDEs driven by non-Gaussian Levy noises
\begin{equation}
dX_t   =  b(X_t )d t + \s(X_t )dL(t),
\end{equation}
have attracted much attention  \cite{Apple, JW, Schertzer} but
this research subject is less developed. Recently, mean exit time
estimates have been investigated by Imkeller et al
\cite{ImkellerP-06, ImkellerP-08} and Yang and Duan
\cite{YangDuan}.

Further progresses in SDEs driven by non-Gaussian   noises will
benefit the research in predictability in weather and climate
systems with non-Gaussian (which is more common) model
uncertainty.

%%%%%%%%%%%%%%%%%%%%%%%%%%%%%%%%%%%
%%%%%%%%%%%%%%%%%%%%%%%%%%%%%%%%%%%
\section{Systems driven by colored noise} \label{colornoise}

Colored noise, or noise with non-zero correlation in time, has
been considered or used in the physical community \cite{Hanggi3,
Gardiner}. A good candidate for modeling colored noise is the
fractional Brownian motion. A fractional Brownian motion (fBm)
process $B^{H}$, where $H\in(0,1)$ is fixed,  is still a Gaussian
process. But it is characterized by the stationarity of its
increments and a  memory property.  The increments of the
fractional Brownian motion are not independent, except in the
standard Brownian case ($H=\frac12$). Thus it is not a  Markov
process except when $H=\frac12$. Specifically, $B^{H}\left(
0\right) =0$ and $Var\left[ B^{H}\left( t\right) -B^{H}\left(
s\right)  \right] =\left\vert t-s\right\vert ^{2H}$. It also
exhibits power scaling and path regularity properties with
H\"older parameter $H$, which are very distinct from Brownian
motion. The standard Brownian motion is a special fBm with
$H=1/2$.

The stochastic calculus involving   fBm is currently being
developed; see e.g. \cite{Nualart, Viens} and references therein.
This will lead to more advances   in the study of SDEs driven by
colored fBm noise:
\begin{equation}
dX_t   =  b(X_t )d t + \s(X_t )dB^{H}(t).
\end{equation}
Since the fBM  $B^{H}(t)$ is not Markov, the  solution process
$X_t$ is not Markov either. Thus the usual techniques from Markov
processes will not be applicable to the study of SDEs driven by
fBms. However, the random dynamical systems approach, as described
in \S \ref{s2} above, looks promising \cite{Maslowski}. The theory
of RDS, developed by Arnold and coworkers  \cite{Arnold},
describes the qualitative behavior of systems of  stochastic
differential equations in terms of stability, Lyapunov exponents,
invariant manifolds, and attractors.

Further progresses in SDEs driven by colored  noises will benefit
the research in predictability in weather and climate systems with
more general (non-white noise) model uncertainty.

\bigskip

%%%%%%%%%%%%%%%%%%%%%%%%
%%%%%%%%%%%%%%%%%%%%%%%%
%%%%%%     Acknowledgement
%%%%%%%%%%%%%%%%%%%%%%%%%%%%%
{\bf Acknowledgements.}

This work was partly supported by the NSF Grants 0542450 and
0620539.

%%%%%%%%%%%%%%
%%%%%%%%%%%%%%
%%%%%%%%%%%%%%


\begin{thebibliography}{50}

\bibitem{Marsden} R. Abraham, J. E. Marsden and T. Ratiu.
\emph{Manifolds, Tensor Analysis, and Applications}. Second Ed.,
Springer-verlag, New York, 1988.


\bibitem{Allen} M. Allen,D. Frame, J. Kettleborough and
 D. Stainforth,
Model error in weather and climate forecasting.
\emph{Predictability of weather and climate} (eds. Palmer, T. and
Hagedorn, R.), pp. 391–428, Cambridge, UK: Cambridge University
Press, 2006.



\bibitem{Apple}
D. Applebaum.
\newblock {\em {L\'evy Processes and Stochastic Calculus}}.
\newblock {Cambridge University Press},  Cambridge, UK, 2004.



\bibitem{Arnold}
L. Arnold.  {\em Random Dynamical Systems.}
  Springer-Verlag, New York, 1998.


\bibitem{Arnold74} L. Arnold,
{\em Stochastic DifferentiL Equations}, John Wiley \& Sons, New
York, 1974.




\bibitem{Arn00}
L.~Arnold.,
\newblock Hasselmann's program visited: The analysis of stochasticity in
  deterministic climate models.
\newblock In J.-S. von~Storch and P.~Imkeller, editors, {\em Stochastic climate
  models}. pages 141--158, Boston, 2001. Birkh{\"a}user.


\bibitem{Aubin}
J.-P. Aubin and G. Da Prato.
\newblock Stochastic Viability and Invariance.
\newblock {\em  Scuola Norm. Sup. Pisa}, {\bf l27} (1990), 595-694.



\bibitem{Berloff} P. S. Berloff, Random-forcing model of the
mesoscale oceanic eddies. \emph{J. Fluid Mech.} \textbf{529}
(2005), 71-95.

%\bibitem{Berloff2} P.S. Berloff, W. Dewar, S. Kravtsov, J. McWilliams and M.
%Ghil. Ocean eddy dynamics in a coupled ocean-atmosphere model:
%Diagnostics and parameterization. \emph{Preprint}, 2005.


\bibitem{Berselli} L.C. Berselli,  T. Iliescu  and W. J. Layton.
    \emph{Mathematics of Large Eddy Simulation of Turbulent
    Flows}.   Springer Verlag, 2005.


\bibitem{BlomkerDuan} D. Blomker  and J. Duan, Predictability of the Burgers dynamics
under model uncertainty. In Boris Rozovsky 60th birthday volume
\emph{Stochastic Differential Equations: Theory and Applications},
P. Baxendale and S. Lototsky (Eds.), p.71-90, World Scientific,
New Jersey, 2007.


\bibitem{Boffetta} G. Boffetta,   A. Celani, M. Cencini,
 G. Lacorata and A. Vulpiani,
The predictability problem in systems with an uncertainty in the
evolution law.  \emph{J. Phys. A}  \textbf{33} (2000), 1313-1324.


\bibitem{Boffetta2}  G. Boffetta,M. Cencini, M. Falcioni and A.
Vulpiani, Predictability: a way to characterize complexity.
\emph{Phys. Rep.}  \textbf{356} (2002),  no. 6, 367--474.




%\bibitem{Bong} V. P. Bongolan-Walsh, J. Duan, and T. Ozgokmen. Dynamics of
%Transport under Random Fluxes on the Boundary.
%\emph{Communications in Nonlinear Science and Numerical
%Simulation} \textbf{13} (2008), 1627-1641.

\bibitem{BrannanDuanErvin}  J. Brannan, J. Duan and V. Ervin,
Escape Probability, Mean Residence Time and Geophysical Fluid
Particle Dynamics, {\em Physica D} {\bf 133} (1999), 23-33.


\bibitem{CL}
T. Caraballo and J. Langa.
\newblock A comparison of the longtime
behavior of linear Ito and Stratonovich partial differential
equations.
\newblock {\em Stochastic Anal. Appl.}, {\bf 19}(2001), no.2 183-195.



\bibitem{Cessi} P. Cessi and S. Louazel,
Decadal oceanic response to stochastic wind forcing, {\em J. Phys.
Oceanography} {\bf 31} (2001), 3020-3029.






\bibitem{Chorin} A. Chorin, A. Kast and R. Kupferman.
Unresolved computation and optimal predictions.  \emph{Comm. Pure
Appl. Math.}, \textbf{52} (1999), pp. 1231-1254.






\bibitem{Crauel} H. Crauel and M. Gundlach (Eds.), \emph{Stochastic Dynamics}.
Papers from the Conference on Random Dynamical Systems held in
Bremen, April 28--May 2, 1997.   Springer-Verlag, New York, 1999.









\bibitem{DaPrato} G. Da Prato and J. Zabczyk, {\em Stochastic Equations
in Infinite Dimensions}, Cambridge University Press, 1992.



\bibitem{DelSole} T. DelSole. Stochastic models of
quasigeostrophic turbulence. \emph{Surveys in Geophysics}
\textbf{25} (2004), 107-149.

%\bibitem{DelSole-Farrell} T. DelSole and B. F. Farrell,
%A stochastically excited linear system as a model for
%quasigeostrophic turbulence: Analytic results for one- and
%two-layer fluids, {\em J. Atmos. Sci.} {\bf 52} (1995) 2531-2547.



\bibitem{Dit} P. D. Ditlevsen, Observation of $\alpha-$stable
noise induced millennial climate changes from an ice record.
\emph{Geophys. Res. Lett. } \textbf{26} (1999), 1441-1444.



\bibitem{DuDuan1} A. Du and J. Duan, Invariant manifold reduction
  for stochastic  dynamical systems.
  \emph{Dynamical Systems and Applications} \textbf{16}(2007), 681-696.



\bibitem{DuDuan} A. Du and J. Duan. A stochastic approach for parameterizing
  unresolved scales in a system with memory.
  \emph{Submitted}, 2007.



\bibitem{DuanNadiga} J. Duan  and B. Nadiga. Stochastic parameterization for
large eddy simulation of geophysical flows. \emph{Proc. Amer.
Math. Soc.} \textbf{135} (2007), 1187-1196.





%\bibitem{DuanGaoSchm}    J. Duan, H. Gao and B. Schmalfuss,
%Stochastic Dynamics of a Coupled Atmosphere-Ocean Model, {\em
%Stochastics and Dynamics} {\bf 2} (2002), 357--380.



\bibitem{DuanLuSchm}
J. Duan, K. Lu  and B. Schmalfu{\ss}. Invariant manifolds for
stochastic partial differential equations.
  {\em The Annals of Probability}, {\bf 31}(2003), 2109-2135.


\bibitem{DuanLuSchm2} J. Duan, K. Lu and B.  Schmalfu{\ss}.
Smooth stable and unstable manifolds for stochastic evolutionary
equations. {\em J. Dynamics and Diff. Eqns.} {\bf 16} (2004),
949-972.




%\bibitem{Ebeling} W. Ebeling and I. M. Sokolov,
%\emph{Statistical Thermodynamics and Stochastic Theory of
%Nonlinear Systems Far from Equilibrium}. Worldscientific, New
%Jersey, USA, 2005.

\bibitem{Farrell} B. F.  Farrell and P. J. Ioannou, Optimal perturbation of uncertain systems. Special issue on
stochastic climate models. \emph{ Stoch. Dyn.}  \emph{2} (2002),
no. 3, 395--402.


\bibitem{Farrell2} B. F.  Farrell and P. J. Ioannou,   Perturbation growth and
structure in uncertain flows. I, II.  \emph{J. Atmospheric Sci.}
\textbf{59} (2002),  no. 18, 2629--2646, 2647--2664.






\bibitem{Filipovic}
\newblock D. Filipovic.
\newblock Invariant manifolds for weak solutions to  stochastic equations.
\newblock {\em Probability Theory \& Related Fields }, Volume \textbf{118} (2000),
Number 3. 323 - 341.






 \bibitem{Freidlin} M. I. Freidlin and A. D. Wentzell,
 {\em Random Perturbations of Dynamical Systems}. 2nd Edition,
 Springer-Verlag, 1998.







\bibitem{Gar} J. Garcia-Ojalvo and J. M. Sancho,
{\em Noise in Spatially Extended Systems}. Springer-Verlag, 1999.

\bibitem{Gardiner} C. W. Gardiner, \emph{Handbook of Stochastic Methods}.
Second Ed., Springer, New York, 1985.




%\bibitem{Greven} A. Greven, G. Keller and G.  Warnecke (Editors),
%\emph{Entropy}. Princeton University Press, 2003.




 \bibitem{Griffa} A. Griffa and S. Castellari,
 Nonlinear general circulation of an ocean model driven by wind with
 a stochastic component. {\em J. Marine Research}, {\bf 49} (1991),
 53-73.










\bibitem{GH}
J. Guckenheimer and P. Holmes.
\newblock \emph{Nonlinear Oscillations,Dynamical Systems and Bifurcations
of Vector Fields}.
\newblock Springer-Verlag, New York, 1983.





%\bibitem{Hanggi2} P. Hanggi and F. Marchesoni (Eds.),
%Stochastic Systems: From Randomness to Complexity.  Special Issue:
%\emph{Physica A},  \textbf{325}(2003), Nos. 1-2: 1-296.



\bibitem{Hanggi3} P. Hanggi and P. Jung, Colored Noise in Dynamical
Systems. \emph{Advances in Chem. Phys}.,
\textbf{89}(1995),239-326.



 \bibitem{Hasselmann} K. Hasselmann, Stochastic climate models:
 Part I. Theory. {\em Tellus}, {\bf 28} (1976), 473-485.


\bibitem{Holloway} G. Holloway, Ocean circulation: Flow in probability under
statistical dynamical forcing, In {\em Stochastic Models in
Geosystems}, S. Molchanov and W. Woyczynski (eds.),
Springer-Verlag, New York, 1996.


  \bibitem{Horst} W. Horsthemke and R. Lefever,
 {\em Noise-Induced Transitions}, Springer-Verlag, Berlin, 1984.



\bibitem{Stuart} W. Huisinga, C. Schutte and A.M. Stuart,
Extracting macroscopic stochastic dynamics: Model problems.
 \emph{Comm. Pure Appl. Math.}, \textbf{56}2003, 234-269.

\bibitem{Imkeller3} P. Imkeller and A. Monahan (eds.),
\emph{Conceptual Stochastic Climate
Models}. Special Issue:\emph{Stochastics and Dynamics},
 \textbf{2}(2002),no.3.



\bibitem{ImkellerP-06} P. Imkeller and I. Pavlyukevich,
First exit time of SDEs driven by stable L\'evy processes.
\emph{Stoch. Proc. Appl.} \textbf{116} (2006), 611-642.


\bibitem{ImkellerP-08}
P. Imkeller, I. Pavlyukevich and T. Wetzel,   First exit times for
L\'evy-driven diffusions with exponentially light jumps.
arXiv:0711.0982.




\bibitem{JW} A. Janicki and A. Weron,
{\em Simulation and Chaotic Behavior of $\alpha-$Stable Stochastic
Processes}, Marcel Dekker, Inc., 1994.



\bibitem{Kantz} W. Just, H. Kantz, C. Rodenbeck and M. Helm,
Stochastic modelling: replacing fast degrees of freedom by noise.
{\em J. Phys. A: Math. Gen.}, {\bf 34} (2001),3199--3213.



\bibitem{KS} I. Karatzas and S. E. Shreve,
\emph{Brownian Motion and Stochastic Calculus}. Second Ed.,
Springer, New York, 1991.


\bibitem{Hansen} V. M. Khade and J. A. Hansen,
State dependent predictability: Impact of uncertainty dynamics,
uncertainty structure and model inadequacies. \emph{Nonlinear
Processes in Geophys}.  \textbf{11} (2004), 351-362.


\bibitem{Klebaner}  F. C. Klebaner,
 \emph{Introduction to  Stochastic  Calculus with Applications}.
 Imperial College Press, London, 2005.


\bibitem{Kleeman}   R. Kleeman. Measuring dynamical prediction utility using relative
entropy.\emph{ J. Atmos Sci}, \textbf{59}:2057-2072, 2002.






 \bibitem{Kloeden} P. E. Kloeden and E. Platen,
 {\em Numerical solution of stochastic differential equations},
 Springer-Verlag, 1992; second corrected printing 1995.




\bibitem{Kunita}
H. Kunita.
 {\em Stochastic flows and stochastic differential
equations.}
  Cambridge University Press, 1990.


%\bibitem{Lasota} A. Lasota and M. C. Mackey.
%\emph{Chaos, Fratals and Noise --- Stochastic Aspects of
%Dynamics}. Second Edition, Springer-Verlag, New York, 1994.





  \bibitem{Leith} C. E. Leith, Climate response and fluctuation dissipation,
 {\em J. Atmos. Sci.}, {\bf 32} (1975), 2022-2025.


\bibitem{Lin} C. C. Lin and L. A. Segel,
{\em Mathematics Applied to Deterministic Problems in the Natural
Sciences},  SIAM, Philadelphia, 1988.


 \bibitem{Majda-Kleeman} A. Majda, R. Kleeman and D. Cai,   A mathematical
framework for quantifying predictability through relative entropy.
Special issue dedicated to Daniel W. Stroock and Srinivasa S. R.
Varadhan on the occasion of their 60th birthday.  \emph{Methods
Appl. Anal.}  \textbf{9}  (2002),  no. 3, 425--444.



\bibitem{Mao} X. Mao. {\em Stochastic differenntial equations \&
applications.} Horwood Publishing, England, 1997.


 \bibitem{Maslowski} B. Maslowski and B. Schmalfuss,
 Random dynamical systems and stationary solutions of differential equations driven by the fractional Brownian motion.
\emph{Stochastic Anal. Appl.}  \textbf{22}  (2004),  no. 6,
1577--1607.



\bibitem{McOwen}
R.C McOwen.
\newblock \emph{Partial Differential Equantions}.
\newblock Pearson Education, New Jersy, 2003.




\bibitem{Millet} A. Millet and P.-L. Morien,   On implicit and explicit
discretization schemes for parabolic SPDEs in any dimension.
\emph{Stochastic Process. Appl.}  \textbf{115}  (2005),  no. 7,
1073--1106.

\bibitem{MZZ} S.-E. A. Mohammed, T. Zhang and H, Zhao, The
stable manifold theorem for semilinear stochastic evolution
equations and stochastic partial differential equations, Memoirs
of the American Mathematical Society, Vol. 196 (2008), pages
1-105.

\bibitem{Mu1} M. Mu, W. S. Duan and B. Wang,   Conditional nonlinear optimal perturbation
and its applications.  \emph{Nonlinear Processes in Geophys. }
\textbf{10} (2003), 493-501.

\bibitem{Mu2} M. Mu, W. S. Duan and J. Chou,
Recent advances in predictability studies in China (1999-2002).
\emph{Adv. Atmos. Sci. } \textbf{21} (2004), 437-443.



\bibitem{Muller} P. M\"uller, Stochastic forcing of quasi-geostrophic eddies,
in {\em Stochastic Modelling in Physical Oceanography}, R. J.
Adler, P. M\"uller and B. Rozovskii (eds.), Birkh\"auser, 1996.


\bibitem{Naeh} T. Naeh, M. M. Klosek, B. J. Matkowsky and Z. Schuss,
A direct approach to the exit problem, {\em SIAM J. Appl. Math.}
{\bf 50} (1990), 595-627.

\bibitem{Nicolis}  C. Nicolis, Dynamics of Model Error:
The Role of Unresolved Scales Revisited. \emph{Journal of the
Atmospheric Sciences} \textbf{61} (2004), 1740–1753.

\bibitem{NorCah81}
G.~R. North and R.~F. Cahalan.
\newblock Predictability in a solvable stochastic climate model.
\newblock {\em J. Atmos. Sci.}, 38:504--513, 1981.



 \bibitem{Nualart} D. Nualart, Stochastic calculus with respect to the fractional
Brownian motion and applications. \emph{Contemporary Mathematics }
\textbf{336}, 3-39, 2003.


\bibitem{Oksendal}
B. Oksendal.
  {\em Stochastic Differenntial Equations.} Sixth Ed.,
  Springer-Verlag, New York, 2003.


\bibitem{Orrell} D. Orrell, L. Smith, J. Barkmeijer and T. N.
Palmer, Model error in weather forecasting. \emph{Nonlinear
Processes in Geophys.} \textbf{8} (2001), 357-371.

\bibitem{Orrell2} D. Orrell, Role of the metric in forecast error
growth: How chaotic is the weather? \emph{Tellus} \textbf{54A}
(2002), 350-362.



\bibitem{Palmer2006} T. N. Palmer and R. Hagedorn (eds.),
\emph{Predictability of weather and climate}, Cambridge, UK:
Cambridge University Press, 2006.

\bibitem{Palmer2005} T. N. Palmer, G. J. Shutts,
R. Hagedorn, F. J. Doblas-Reyes, T. Jung and M. Leutbecher.
Representing model uncertainty in weather and climate prediction.
\emph{Annu. Rev. Earth Planet. Sci.}  \textbf{33} (2005), 163-193.


\bibitem{Palmer} T. N. Palmer. A nonlinear dynamical perspective
on model error: A proposal for non-local stochastic-dynamic
parameterization in weather and climate prediction models. {\em Q.
J. Meteorological Soc.} \textbf{127} (2001) Part B, 279-304.



\bibitem{Pas} C. Pasquero and E. Tziperman,  Statistical parameterization
of heterogeneous oceanic convection, \emph{J. Phys. Oceanography},
\textbf{37} (2007), 214-229.





\bibitem{PeiOor92}J.~P. Peixoto and A.~H. Oort, {\em Physics of Climate}.
 Springer, New York, 1992.

\bibitem{Perko}
 L. Perko.
 \newblock {\em Differential Equations and Dynamical
Systems.}
\newblock Cambridge University Press, 1990.



\bibitem{Risken} H. Risken, {\em The Fokker-Planck Equation},
Springer-Verlag, New York, 1984.

\bibitem{Roberts} A. J.  Roberts,   A step towards holistic discretisation of
stochastic partial differential equations.  \emph{ANZIAM J.}
\textbf{45} (2003/04),  (C), C1--C15.


 \bibitem{Roz} B. L. Rozovskii,
{\em Stochastic Evolution Equations}. Kluwer Academic Publishers,
Boston, 1990.


\bibitem{Samelson} R. M. Samelson, Stochastically forced current
fluctuations in vertical shear and over topography, {\em J.
Geophys. Res.} {\bf 94} (1989), 8207-8215.

%\bibitem{Saravanan-McWilliams2} R. Saravanan and J. C. McWilliams,
%Multiple equilibria, natural  variability, and climate transitions
%in an idealized ocean-atmosphere model, {\em J. Climate} {\bf 8}
%(1995), 2296-2323.

\bibitem{Saravanan-McWilliams} R. Saravanan and J. C. McWilliams,
Advective ocean-atmosphere interaction: An analytical stochastic
model with implications for decadal variability, {\em J. Climate}
{\bf 11} (1998), 165-188.


\bibitem{Sato-99}
K.-I. Sato.
\newblock {\em {L{\'e}vy Processes and Infinitely Divisible Distributions}},
\newblock {Cambridge University Press}, Cambridge, 1999.



\bibitem{Schneider} T. Schneider and S. M. Griffies,
A conceptual framework for predictability studies. J. of Climate
\textbf{12} (1999), 3133-3155.

\bibitem{Schertzer}
D. Schertzer, M. Larcheveque, J. Duan,  V. Yanovsky and S.
Lovejoy, Fractional Fokker--Planck Equation for  Nonlinear
Stochastic Differential Equations Driven by Non-Gaussian Levy
Stable Noises. {\em  J. Math. Phys.}, \textbf{42} (2001), 200-212.

\bibitem{Schmalfuss} B. Schmalfuss,
The random attractor of the stochastic Lorenz system.
\emph{Zeitschrift fürangewandte mathematik und physik (ZAMP)}
\textbf{48} (1997), 951 - 975.


\bibitem{Schuss} Z. Schuss, {\em Theory and Applications of
Stochastic Differential Equations}, Wiley $\&$ Sons, New York,
1980.

\bibitem{Shlesinger} M. F. Shlesinger, G. M. Zaslavsky  and U.
Frisch, {L\'evy Flights and Related Topics in Physics} (Lecture
Notes in Physics, 450. Springer-Verlag, Berlin, 1995).


\bibitem{Smith} L. A. Smith, C. Ziehmann and K. Fraedrich,
Uncertainty dynamics and predictability in chaotic systems.
\emph{Q. J. R. Meteorol. Soc.} \textbf{125} (1999), 2855-2886.





%\bibitem{Sura2} P. Sura, K. Fraedrich and F. Lunkeit,
%Regime transitions in a stochastically forced Double-gyre model.
 %{\em J. Phys. Oceanography},\textbf{31}(2001), 411-426.

\bibitem{Sura3} P. Sura  and C. Penland,  Sensitivity of a
double-gyre model to details of stochastic forcing. \emph{Ocean
Modelling}, \textbf{4}(2002), 327-345.


\bibitem{Viens}  C.A. Tudor and F. Viens, Statistical aspects of the
fractional stochastic calculus.\emph{ Annals of Statistics}, Vol.
\textbf{35} (3) (2007), 1183-1212.



%\bibitem{VanKampen2} N. G. Van Kampen,
%How do stochastic processes enter into physics? {\em Lecture Note
%in Phys.} {\bf 1250} (1987)  128--137.


\bibitem{Toth}  S. Vannitsem and Z. Toth,
Short-term dynamics of model errors. \emph{ J. Atmospheric Sci.}
\textbf{59}  (2002),  no. 17, 2594--2604.





\bibitem{WangDuan} W. Wang and J. Duan, A dynamical approximation  for
stochastic partial differential equations,
  \emph{J. Math. Phys.}   \textbf{48}(2007), No. 10, 102701.

\bibitem{Wanner2} T. Wanner.
\newblock  Linearization of Random Dynamical Systems.
\newblock {\em Dynamics Report}, Volumn 4. Spring-Verlog, New York, 1995.

\bibitem{WaymireDuan} E. Waymire  and J. Duan (Eds.).
\emph{Probability and Partial Differential Equations in Modern
Applied Mathematics}. Springer-Verlag,   2005.







%\bibitem{Wiggins} S. Wiggins.
%\newblock \emph{Introduction to Applied Nonlinear Dynamical Systems and Chaos}.
%\newblock Spring-Verlag, New York, 1990.

\bibitem{Wilks} D. S. Wilks, Effects of stochastic
parameterizations in the Lorenz '96 system. \emph{Q. J. R.
Meteorol. Soc.} \textbf{131} (2005), 389-407.


\bibitem{Williams} P. D. Williams. Modelling climate change: the
role of unresolved processes. \emph{Phil. Trans. R. Soc. A} (2005)
\textbf{363}, 2931-2946.

\bibitem{YangDuan} Z. Yang and J. Duan,
Mean exit time estimates for dynamical systems with non-Gaussian
L\'evy noises. \emph{Preprint}, 2008.


\bibitem{Zhao} H. Zhao and Z. Zheng,
Random periodic solutions of random dynamical systems.
\emph{Preprint}, 2008.


\end{thebibliography}
\end{document}